\documentclass[10pt]{article}
\usepackage{amsthm,amsfonts,amsmath,amscd,amssymb,epsfig,verbatim,
}
\usepackage[all]{xy}
\let\oldmarginpar\marginpar
\renewcommand\marginpar[1]{\-\oldmarginpar[\raggedleft\footnotesize #1]%
{\raggedright\footnotesize #1}}

\theoremstyle{plain}
\newtheorem{theorem}{Theorem}
\newtheorem{thm}[theorem]{Theorem}

\newtheorem{cor}[theorem]{Corollary}

\newtheorem{lemma}[theorem]{Lemma}

\theoremstyle{remark}
\newtheorem{conj}[theorem]{Conjecture}

\newtheorem{remark}[theorem]{Remark}
\newtheorem*{remark*}{Remark}
\newtheorem{example}[theorem]{Example}
\newtheorem*{example*}{Example}

\theoremstyle{definition}

\newcommand{\id}{\mathrm{Id}}

\newcommand{\wt}{\widetilde}
\newcommand{\wh}{\widehat}
\newcommand{\ol}{\overline}

\newcommand{\p}{\partial}

\newcommand{\om}{\omega}

\newcommand{\eps}{\varepsilon}

\newcommand{\la}{\langle}
\newcommand{\ra}{\rangle}

\newcommand{\R}{{\mathbb{R}}}
\newcommand{\C}{{\mathbb{C}}}

\newcommand{\Crit}{{\rm Crit}}

\newcommand{\Zero}{\mathrm{Zero}}

\newcommand{\Stein}{\mathfrak{Stein}}
\newcommand{\Weinstein}{\mathfrak{Weinstein}}

\newcommand{\Functions}{\mathfrak{Functions}}

\newcommand{\FF}{\mathcal{F}}

\newcommand{\GG}{\mathcal{G}}

\newcommand{\WW}{\mathcal{W}}

\newcommand{\fS}{\mathfrak{S}}

\newcommand{\fW}{{\mathfrak W}}

\newcommand{\grad}{\rm grad}
\title{A note on gradient-like vector fields}
\author{Kai Cieliebak}
\date{}

\begin{document}
\maketitle

\begin{abstract}
This note proposes a new notion of a gradient-like vector field and
discusses its implications for the theory of Stein and Weinstein structures.
\end{abstract}

For a smooth function $\phi:V\to\R$ on a manifold $V$, its gradient
vector field $\grad_g\phi$ with respect to a Riemannian metric $g$ is defined
by $d\phi=g(\grad_g\phi,\cdot)$. For many applications, this notion is
too rigid and one wants to allow more general ``gradient-like'' vector
fields $X$ satisfying 
\begin{equation}\label{eq:weak-grad}
  d\phi(X)>0\text{ outside }\Zero(X)=\Crit(\phi),
\end{equation}
where $\Zero(X)$ denotes the zero set of $X$ and $\Crit(\phi)$ the set
of critical points of $\phi$. 
For example, this occurs in Morse theory~\cite{Mil65}, and in the
theory of Weinstein structures where $X$ is additionally required to
be Liouville for a symplectic form~\cite{CE12}.
The drawback of this notion is that the zeroes of $X$ can be of very
different nature from the critical points of $\phi$, so one does not
expect a good deformation theory for pairs $(X,\phi)$ satisfying~\eqref{eq:weak-grad}. 
To remedy this, the following more quantitative notion
of~\eqref{eq:weak-grad} is used e.g.~in~\cite{CE12}:
\begin{equation}\label{eq:Lyap}
  d\phi(X)\geq\delta\cdot (|X|^2+|d\phi|^2)\text{ for a function }\delta:V\to\R_{>0}.
\end{equation}
However, the deformation theory of such pairs is still not
understood. For example, it is unknown whether any pair $(X,\phi)$
satisfying~\eqref{eq:Lyap} can be deformed through such pairs to one
for which the function is Morse. 
Because of this, the map from Weinstein to Stein structures in~\cite{CE12}
is only defined under the condition that the involved functions are
generalized Morse. While this suffices for the discussion of $0$- and
$1$-parametric familes of such structures in~\cite{CE12}, it has been
an obstacle to extending these results to higher dimensional families.

It is the goal of this note to show that these issues get resolved by
the following notion. Let us call a vector field $X$ {\em
  gradient-like} for $\phi$ if
\begin{equation}\label{eq:grad-like}
  d\phi=g(X,\cdot) \text{ for some positive smooth $(2,0)$ tensor field $g$,}
\end{equation}
where ``positive'' means $g(v,v)>0$ for all $v\neq 0$. If $g$ is
in addition symmetric this just means $X=\grad_g\phi$, but in general
symmetry of $g$ is not required. 
Denoting $\mathfrak{X}$ the space of vector fields, $\FF$ the
space of smooth functions $\phi:V\to\R$, and $\GG$ the space of
positive $(2,0)$ tensor fields, we have the commuting diagram 
\begin{equation*}
\xymatrix
@C=0pt{
  \mathfrak{X} & \ar[l]_{\rm contr\ \ \ \ \ \ \ \ \ \ \ \ \ \ \ \ \ \ \ \ \ \ }
  \{(X,\phi,g)\in\mathfrak{X}\times\FF\times\GG\mid
  d\phi=g(X,\cdot)\} \ar[dl]^{\rm contr} \ar[dr]_\cong \\
  \{(X,\phi)\in\mathfrak{X}\times\FF \mid \eqref{eq:grad-like}\}
  \ar[u] \ar[dr]_{\rm Serre} & & \FF\times\GG \ar[dl]^{\rm proj} \\
  & \FF
}
\end{equation*}
Here the middle right arrow is a homeomorphism (because
$d\phi=g(X,\cdot)$ uniquely determines $X$ given $\phi,g$),  
the upper and middle arrow have contractible fibres, and the lower right
arrow is the canonical projection $(\phi,g)\mapsto\phi$. So we obtain

\begin{lemma}
The lower left arrow in the above diagram is a Serre fibration with
contractible fibres and a 
section (just fixing some $g$), so every family of functions lifts to
a family of pairs satisfying~\eqref{eq:grad-like}. In
particular, any pair $(X,\phi)$ satisfying~\eqref{eq:grad-like} is
homotopic through such pairs to one for which the function is Morse,
and any $1$-parametric family of such pairs is homotopic to one having
only birth-death critical points. \qed
\end{lemma}

Recall that a critical point $p$ of $\phi$ is called {\em Morse} if
the Hessian of $\phi$ at $p$ is nondegenerate, and {\em embryonic} if
the Hessian has a $1$-dimensional kernel on which the third derivative
of $\phi$ does not vanish. The following lemma shows that at such critical
points the notions~\eqref{eq:Lyap} and~\eqref{eq:grad-like} are
equivalent. 

\begin{lemma}\label{lem:Morse-embryonic}
(a) If $(X,\phi)$ satisfies~\eqref{eq:grad-like}, then it also
  satisfies~\eqref{eq:Lyap}. \\
(b) If $(X,\phi)$ satisfies~\eqref{eq:Lyap} and $\phi$ has only Morse or
embryonic critical points, then $(X,\phi)$ also
satisfies~\eqref{eq:grad-like}. \\
(c) If $(X,\phi)$ satisfies~\eqref{eq:grad-like} near $\Crit(\phi)$
and~\eqref{eq:weak-grad}, then it satisfies~\eqref{eq:grad-like}
everywere.   
\end{lemma}

\begin{proof}
(a) Given a Riemannian metric with norm $|\ |$, define two functions $a,B:V\to(0,\infty)$ by
$$
  a(p):=\inf_{0\neq v\in T_pV}\frac{g_p(v,v)}{|v|^2},\qquad
  B(p):=\sum_{0\neq v,w\in T_pV}\frac{g_p(v,w)}{|v||w|}.
$$
Then $d_p\phi(X)=g_p(X,X)\geq a(p)|X|^2$ and $|d_p\phi(v)|=|g_p(X,v)|\leq
B(p)|X||v|$ combine to the following inequalities which
imply~\eqref{eq:Lyap}: 
$$
  |d_p\phi|^2 = \sum_{0\neq v\in T_pV}\frac{|d_p\phi(v)|^2}{|v|^2}
  \leq B(p)^2|X|^2\leq \frac{B(p)^2}{a(p)^2}d_p\phi(X).
$$
(b) {\em Step 1. }Near a Morse singularity, by the Morse lemma we can choose local
coordinates $Z\in U\subset\R^m$ in which 
$$
   \phi(Z)=\phi(0)+\frac12\la BZ,Z\ra\quad\text{and}\quad
   X(Z) = A(Z)Z
$$
for a symmetric invertible matrix $B$ and a smooth map
$A:U\to\R^{m\times m}$ with $A(0)=A_0$ invertible.
Condition~\eqref{eq:Lyap} implies that
$\la Bv,A_0v\ra\geq\beta|v|^2$ for some $\beta>0$.
We are looking for a $(2,0)$-tensor field $g_Z(v,w)=\la G(Z)v,w\ra$ satisfying
$$
   \la G(Z)X(Z),\cdot\ra = d\phi(Z) = \la BZ,\cdot\ra,
$$
or equivalently, $G(Z)X(Z)=G(Z)A(Z)Z=BZ$. This is
satisfied if we define $G:U\to\R^{m\times m}$ by
$$
   G(Z)w := BA(Z)^{-1}w.
$$
Writing $w=A(Z)v$, for $U$ sufficiently small we then have
$$
   \la G(Z)w,w\ra = \la Bv,A(Z)v\ra\geq\frac{\beta}{2}|v|^2\geq\gamma|w|^2
$$
for some $\gamma>0$, so $g$ is positive. 

{\em Step 2. }Near an embryonic singularity, by Theorem 9.4 and the proof of
Lemma 9.12 in~\cite{CE12} we can choose local coordinates $Z=(w,z)\in
U\subset\R^{m-1}\times\R$ in which 
$$
   \phi(Z)=\phi(0)+\frac12\la Bw,w\ra + \frac13 cz^3\quad\text{and}\quad
   X(Z) = \Bigl(A(Z)w,a_1(Z)z^2 + a_2(Z)w\Bigr)
$$
for a symmetric invertible matrix $B$, a number $c>0$, and smooth maps
$A:U\to\R^{(m-1)\times(m-1)}$ with $A(0)=A_0$ invertible, $a_1:U\to\R$ with
$a_1(0)=1$, and $a_2:U\to\R^{m-1}$ with $a_2(0)=0$. 
Condition~\eqref{eq:Lyap} implies that
$\la Bv,A_0v\ra\geq\beta|v|^2$ for some $\beta>0$.
We are looking for a $(2,0)$-tensor field $g_Z(v,w)=\la G(Z)v,w\ra$ satisfying
$$
   \la G(Z)X(Z),\cdot\ra = d\phi(Z) = \la (Bw,cz^2),\cdot\ra.
$$
We will look for $G$ in the form
$$
G(Z) = \begin{pmatrix}
   g_{11}(Z) & 0 \\ g_{21}(Z) & g_{22}(Z)
\end{pmatrix}
$$
with respect to the splitting $\R^m=\R^{m-1}\times\R$. Then the
preceding equation is equivalent to
$$
   g_{11}(Z)A(Z)w = Bw\quad\text{and}\quad
   g_{21}(Z)A(Z)w + g_{22}(Z)\Bigl(a_1(Z)z^2 + a_2(Z)w\Bigr) = cz^2.
$$
We solve the first equation by setting $g_{11}(Z):=BA(Z)^{-1}$. To
solve the second equation, we first set $g_{22}(Z):=ca_1(Z)^{-1}$ and then
$g_{12}(Z):=-g_{22}(Z)a_2(Z)A(Z)^{-1}$. 
By the assumptions $a_2(0)=0$ etc, we can estimate for $U$ sufficiently small:
\begin{align*}
   \la G(Z)\begin{pmatrix}\wh w \\ \wh
     z\end{pmatrix},\begin{pmatrix}\wh w \\ \wh z\end{pmatrix}\ra  
   &= \la g_{11}(Z)\wh w,\wh w\ra + \la g_{21}(Z)\wh w,\wh z\ra
       + \la g_{22}(Z)\wh z,\wh z\ra \cr
   &= \la BA(Z)^{-1}\wh w,\wh w\ra
       - g_{22}(Z)a_2(Z)A(Z)^{-1}\wh w\wh z + ca_1(Z)^{-1}\wh z^2 \\
   &\geq \gamma|\wh w|^2 - C|Z|\,|\wh w|\,|\wh z| + \frac{c}{2}\wh z^2
\end{align*}
for some constants $\gamma,c,C>0$. This shows that $g$ is positive for
$U$ sufficiently small. 

{\em Step 3. }The preceding two steps provide $g$ in a neighbourhood of the
critical points. Away from the critical points the existence of $g$ is
obvious, and the constructions on the various regions can be patched
together using a partition of unity. 

(c) By hypothesis there exist open sets $U_0,U_1$ with
$\Crit(\phi)\subset U_0\subset\ol U_0\subset U_1$ and a positive $(2,0)$
tensor field $g_1$ on $U_1$ such that $d\phi=g_1(X,\cdot)$ on $U_1$.
Since $d\phi(X)>0$ on $V\setminus U_0$, there exists a Riemannian
metric $g_0$ on $V\setminus U_0$ such that $d\phi=g_0(X,\cdot)$ on
$V\setminus U_0$ (take any $g_0$ for which $X$ is orthogonal to
$\ker d\phi$ and $g_0(X,X)=d\phi(X)$). Now pick a cutoff function
$\chi:V\to[0,1]$ which equals $1$ on $U_0$ and $0$ on $V\setminus U_1$
and define $g$ by $g_1$ on $U_0$, $g_0$ on $V\setminus U_1$, and
$(1-\chi)g_0+\chi g_1$ on $U_1\setminus U_0$. 
\end{proof}

\begin{remark}
For a pair $(X,\phi)$ we have the implications
$(4) \Longrightarrow (3) \Longrightarrow (2) \Longrightarrow (1)$
with
\begin{equation}\label{eq:grad}
  X=\grad_g\phi \text{ for some Riemannian metric $g$,}
\end{equation}
and the converse implications do not hold. For example, on $\R$ the pair
$\phi(x)=x^3$ and $X(x)=x^4$ satisfies (1) but not (2), and the pair
$$
\phi(x) = \begin{cases}
  e^{-\frac{1}{x^2}} & x\neq 0, \cr
  0 & x=0,
\end{cases}
\qquad
X(x) = \begin{cases}
  \frac{2}{x^3}e^{-\frac{1}{x^2}} & x< 0, \cr
  0 & x=0, \cr
  \frac{1}{x^3}e^{-\frac{1}{x^2}} & x> 0
\end{cases}
$$
satisfies (2) but not (3). A pair $(X,\phi)$ satisfying (3) but not
(4) is given in~\cite[Example 11.18]{CE12}. 
\end{remark}



\begin{example}
In dimension $1$ for $(X,\phi)$ real analytic, (2) implies (3). As
pointed out by Y.~Eliashberg, this is no longer the case in higher
dimensions. For example, the pair $\phi(x,y)=\frac14(x^4+y^4)$ and
$X(x,y)=(x^3+x^2y^2)\p_x+y^3\p_y$ satisfies (2) but not (3). To see
this, note that $X=\nabla\phi+F\p_x$ with the function $F(x,y)=x^2y^2$.
It does not satisfy (3) because $F$ does not lie in the ideal $I$
generated by the partial derivatives $\p_x\phi=x^3$ and $\p_y\phi=y^3$
(over the ring of germs of smooth functions), and it satisfies (2)
because $F^2=xy\p_x\phi\p_y\phi\in I^2$.  
\end{example}

{\bf Application to Weinstein structures. }
Let us now use the new notion of ``gradient-like'' to define a {\em
  Weinstein structure} on $V$ as a triple $(\om,X,\phi)$ where

\ \ \ \ \ (i) $\om$ is a symplectic form;

\ \ \ \ \ (ii) $X$ is a Liouville vector field for $\om$, i.e.~$L_X\om=\om$;

\ \ \ \ \ (iii) $\phi:V\to\R$ is a smooth exhausting function such that
$(X,\phi)$ satisfies~\eqref{eq:grad-like}.

Here ``exhausting'' means that $\phi$ is proper and bounded from below. 
The form $\lambda=i_X\om$ satisfies $d\lambda=\om$ and is called the
Liouville form. 

Note that there is no condition on the kinds of critical points of
$\phi$. In contrast, the definition of a ``Weinstein structure''
in~\cite{CE12} used instead of (iii) the condition that $\phi$ has only
Morse or embryonic critical points and $(X,\phi)$ satisfies~\eqref{eq:Lyap}.
By Lemma~\ref{lem:Morse-embryonic}, $(X,\phi)$ then also
satisfies~\eqref{eq:grad-like}, so the above definition of a Weinstein
structure generalizes the one in~\cite{CE12}. 

The above definition of a Weinstein structure includes all the
natural examples known to me:

\begin{example}\label{ex:Weinstein}
(a) Recall that a {\em Stein structure} $(J,\phi)$ consists of a
complex structure $J$ and a smooth exhausting function $\phi:V\to\R$
which is $J$-convex, i.e.~$-d(d\phi\circ J)(v,Jv)>0$ for all $v\neq 0$.
It gives rise to a natural Weinstein structure $\WW(J,\phi)=(\om_\phi,X_\phi,\phi)$,
where $\om_\phi=-d(d\phi\circ J)$ and $X_\phi$ is the gradient of
$\phi$ with respect to the Riemannian metric $g_\phi=\om_\phi(\cdot,J\cdot)$.

(b) The construction in (a) still works if $J$ is just an almost
complex structure, the only difference being that $g_\phi$ need not be
symmetric, so the pair $(X_\phi,\phi)$ satisfies~\eqref{eq:grad-like}
but not necessarily~\eqref{eq:grad}. See~\cite[Remark 11.17]{CE12}. 

(c) The cotangent bundle $T^*M$ of a closed Riemannian manifold $M$
carries the canonical Weinstein structure $(\om,X,\phi)$ where
$\om=d\lambda$ for the canonical Liouville form $\lambda=p\,dq$, $X$
is defined by $i_X\om=\lambda$, and $\phi(q,p)=\frac12|p|^2$ for the
dual metric. Here $X$ is the gradient of $\phi$ with respect to the
Riemannian metric on $T^*M$ induced by the metric on $M$ and its
Levi-Civita connection.

(d) Weinstein manifolds with {\em arboreal skeleton} also fall into
this class, see~\cite{AEN24}. 
\end{example}

On the other hand, the above definition is restrictive enough to allow
for a good deformation theory:

\begin{lemma}\label{lem:Weinstein}
Let $(V,\om,X,\phi)$ be a Weinstein manifold with Liouville form
$\lambda$. Then for each smooth function 
$\wt\phi$ which is sufficiently close to $\phi$ in the strong
$C^2$-topology there exists a Weinstein homotopy from
$(\om,X,\phi)$ to a Weinstein structure $(\wt\om,\wt X,\wt\phi)$. 
If each connected component of the support $S$ of $\wt\phi-\phi$ is compact, then
we can choose the Weinstein homotopy fixed outside $S$, with fixed $\om_t=\om$, and
with Liouville forms $\lambda_t$ such that $\lambda_t-\lambda$ is exact. 
In particular, $(\om,X,\phi)$ is Weinstein homotopic with fixed $\om$
to a Weinstein structure $(\om,\wt X,\wt\phi)$ with $\wt\phi$ Morse and
$\wt\phi=\phi$ outside a neighbourhood of $\Crit(\phi)$.  
\end{lemma}

\begin{proof}
Let $\lambda=i_X\om$. Define a smooth bundle homomorphism $A:TV\to TV$ by
$\om(\cdot,A\cdot)=g(\cdot,\cdot)$. Then $A$ is invertible and
$$
   d\phi=g(X,\cdot)=\om(X,A\cdot)=\lambda(A\cdot)
$$
shows $\lambda=d\phi\circ A^{-1}$. 

Let now $\wt\phi$ be close to $\phi$ in the strong $C^2$-topology
(see~\cite{Hir76}). 
Then $\wt\lambda:=d\wt\phi\circ A^{-1}$ is $C^1$-close to $\lambda$ and
$\wt\om:=d\wt\lambda$ is $C^0$-close to $\om$, hence symplectic. 
The vector field $\wt X$ defined by $i_{\wt X}\wt\om=\wt\lambda$ is
$C^0$-close to $X$ and satisfies $L_{\wt X}\wt\om=\wt\om$. Now
$$
   d\wt\phi = \wt\lambda(A\cdot) = \wt\om(\wt X,A\cdot) = \wt g(\wt X,\cdot)
$$
with the $(2,0)$-tensor field $\wt g:=\wt\om(\cdot,A\cdot)$. Since
$\wt g$ is $C^0$-close to $g$, it is also positive (if
$g_x(v,v)\geq\delta>0$ for all $v$ with $|v|=1$, then $\wt
g_x(v,v)\geq\delta/2>0$ for all $v$ with $|v|=1$). Hence $(\wt
X,\wt\phi)$ satisfies~\eqref{eq:grad-like}, and therefore $(\wt\om,\wt
X,\wt\phi)$ is a Weinstein structure.  
Applying the same argument to the functions
$\phi_t=(1-t)\phi+t\wt\phi$, $t\in[0,1]$, we find a Weinstein homotopy
$(\om_t,X_t,\phi_t)$ (in the sense of~\cite{CE12}) from $(\om,X,\phi)$
to $(\wt\om,\wt X,\wt\phi)$. 

Suppose now that each connected component of the support $S$ of
$\wt\phi-\phi$ is compact. Then Moser's theorem as stated
in~\cite[Theorem 6.8]{CE12}) yields a family of diffeomorphisms
$h_t:V\to V$ with $h_0=\id$ and $h_t=\id$ outside $S$ such that
$h_t^*\lambda_t-\lambda$ is exact. Thus $(\om,h_t^*X_t,h_t^*\phi_t)$
is the desired Weinstein homotopy. 

Finally, by Sard's theorem there exists a neighbourhood $S$ of
$\Crit(\phi)$ with each connected component compact, so the last
assertion follows by taking $\wt\phi$ to be Morse with $\wt\phi=\phi$
outside $S$.
\end{proof}

This definition of a Weinstein structure provides an appropriate
setting for the results in~\cite{CE12}. 
Fix a compact smooth manifold $W$ with boundary and assume that all
considered functions on $W$ have $\p W$ as their regular maximum level set.
Denote by $\Stein$ and $\Weinstein$ the spaces of Stein and Weinstein
domain structures on $W$, respectively, and by $\Functions$ the space of
all smooth functions $\phi:W\to\R$.
Note that these spaces involve no nondegeneracy assumptions on the
underlying functions, whereas in~\cite{CE12} it was assumed that all
the underlying functions are generalized Morse. 
We have the commutative diagram
$$
   \xymatrix{\Stein     \ar[rd]^{\pi_{\fS} }\ar[rr]^{\fW}& &\Weinstein
     \ar[ld]_{\pi_{\fW}}\\ 
   &\Functions &}
$$
where $\fW(J,\phi)=(\om_\phi,X_\phi,\phi)$ as in
Example~\ref{ex:Weinstein}(a) above,  
$\pi_\fW(\om,X,\phi):=\phi$, and $\pi_\fS(J,\phi):=\phi$. 

Following~\cite{Gro86} let us call a continuous map $p:E\to B$
between topological spaces a {\em microfibration} if, given 
continuous maps $f:D^k\times[0,1]\to B$ and $F_0:D^k\to E$ with
$p\circ F_0=f|_{D^k\times 0}$, there exists $\eps\in(0,1]$ and a
continuous lift $F:D^k\times[0,\eps]\to E$ such that $F|_{D^k\times
  0}=F_0$ and $p\circ F=f|_{D^k\times[0,\eps]}$. 
If we can choose $\eps=1$ it is a {\em Serre fibration}. 

\begin{cor}\label{cor:Weinstein}
The map $\pi_{\fW}:\Weinstein\to\Functions$ is a microfibration.
\end{cor}

\begin{proof}
Let continuous families of Weinstein structures $(\om_p,X_p,\phi_p)$
and functions $\phi_{p,t}$, $(p,t)\in D^k\times[0,1]$, with
$\phi_{p,0}=\phi_0$ be given. Since $W$ is compact, there exists
$\eps\in(0,1]$ such that for each $(p,t)\in D^k\times[0,\eps]$ the
function $\phi_{p,t}$ is sufficiently $C^2$-close to $\phi_p$ so that
Lemma~\ref{lem:Weinstein} applies. Since the construction in its proof is
completely canonical, it provides the desired continuous lift
$(\om_{p,t},X_{p,t},\phi_{p,t})$ for $(p,t)\in D^k\times[0,\eps]$.
\end{proof}

\begin{remark}
(a) The map $\pi_{\fW}:\Weinstein\to\Functions$ is {\em not} a Serre
fibration. For example, for $\phi$ in the image of $\pi_{\fW}$ the
unstable manifold at each critical point must have at least dimension
$\dim W/2$ (see~\cite[Proposition 11.9]{CE12}). 
It would be interesting to know whether, for a function with isolated
critical points, this is the only condition for being in the image of
$\pi_{\fW}$. \\
(b) Lemma~\ref{lem:Weinstein} also yields an analogue of
Corollary~\ref{cor:Weinstein} for Weinstein manifolds instead of
Weinstein domains, which is however more complicated to state.
\end{remark}

In view of Corollary~\ref{cor:Weinstein}, Theorem 1.1 in~\cite{CE12}
becomes

\begin{thm}[\cite{CE12}]
The map $\fW:\Stein\to\Weinstein$ induces an isomorphism on $\pi_0$
and a surjection on $\pi_1$. 
\end{thm} 

Conjecture 1.4 in~\cite{CE12} now takes the form

\begin{conj}[\cite{CE12}]
The map $\fW:\Stein\to\Weinstein$ is a weak homotopy equivalence.   
\end{conj}

{\bf Hypersurfaces in contact manifolds. }
For background on the following discussion see~\cite{HH19}.
Consider a compact cooriented codimension $1$ hypersurface $\Sigma$ in a
cooriented contact manifold $(M,\xi=\ker\alpha)$. It carries a
$1$-dimensional characteristic foliation with singularities at points
$p$ where $\xi_p=T_p\Sigma$, where $p$ is called positive
if the coorientations of $\xi_p$ and $T_p\Sigma$ agree and negative otherwise. 
The restriction $\beta=\alpha|_\Sigma$ is a Liouville form near the
singular points, and we can choose a vector field $X$ generating the
characteristic foliation which agrees with $\pm$ the Liouville field
near $\pm$ singular points. Let us call $\Sigma$ {\em gradient-like}
if $X$ is gradient-like for some function $\phi:\Sigma\to\R$ in the
sense of~\eqref{eq:grad-like}. If $\phi$ has only Morse resp.~Morse
and birth-death critical points this agrees with the notion of a 
Morse resp.~$1$-Morse hypersurface in~\cite{HH19}. 
Since every exact deformation of $\beta$ near its singular locus can
be realized by a deformation of $\Sigma$ in $(M,\alpha)$,
Lemma~\ref{lem:Weinstein} yields 

\begin{cor}\label{cor:hypersurface}
The map $(\Sigma,X,\phi)\mapsto\phi$ that associates to each
gradient-like hypersurface its function is a microfibration. In particular, each
gradient-like hypersurface can be $C^2$-perturbed to one with only
Morse singularities. \qed
\end{cor}

{\bf Acknowledgement. }
The idea of this note arose during preparation for a ``Mathematical
Conversation'' at the Institute for Advanced Study on January 26, 2022. 
I thank its organisers and participants for making this happen despite
freezing temperatures, and Y.~Eliashberg for inspiring conversations
afterwards.



\begin{thebibliography}{9999999}

\bibitem{AEN24}
D.~\'Alvarez-Gavela, Y.~Eliashberg and D.~Nadler, {\em Weinstein
  manifolds as cotangent buildings}, in preparation. 
  
\bibitem{CE12}
K.~Cieliebak and Y.~Eliashberg, {\em From Stein to Weinstein and
  Back. Symplectic Geometry of Affine Complex Manifolds}, American
Mathematical Society Colloquium Publications 59 (2012).  

\bibitem{Gro86}
M.~Gromov, {\em Partial Differential Relations}, Springer (1986).

\bibitem{Hir76} M.~Hirsch, {\em Differential Topology},
Springer (1976).

\bibitem{HH19}
K.~Honda and Y.~Huang, {\em Convex hypersurface theory in contact
  topology}, arXiv:1907.06025.

\bibitem{Mil65}
J.~Milnor, {\em Lectures On The H-Cobordism Theorem}, Notes by
L.~Siebenmann and J.~Sondow, Princeton University Press (1965). 
\end{thebibliography}
\end{document}